\documentclass[10pt, a4paper, twoside,reqno]{amsart}

\usepackage{fancyhdr}

\makeatletter

\def\section{\@startsection{section}{1}%
	\z@{.7\linespacing\@plus\linespacing}{.5\linespacing}%
	{\normalfont \Large\scshape\centering}}

\def\subsection{\@startsection{subsection}{2}%
	\z@{.5\linespacing\@plus.7\linespacing}{.5\linespacing}%
	{\normalfont\large\bfseries}}

\def\subsubsection{\@startsection{subsubsection}{3}%
	\z@{.5\linespacing\@plus.7\linespacing}{.5\linespacing}%
	{\normalfont\itshape}}

\usepackage[usenames, dvipsnames]{color}
\definecolor{darkblue}{rgb}{0.0, 0.0, 0.45}

\usepackage[colorlinks	= true,
raiselinks	= true,
linkcolor	= darkblue, 
citecolor	= Mahogany,
urlcolor	= ForestGreen,
pdfauthor	= {Georgios Darivianakis},
pdftitle	= {},
pdfkeywords	= {},
pdfsubject	= {},
plainpages	= false]{hyperref}

\allowdisplaybreaks
\date{\today}

\usepackage{graphicx, psfrag}
\usepackage{amsthm}
\usepackage{amsmath} 
\usepackage{amssymb}  
\usepackage[noadjust]{cite}
\usepackage{color}
\usepackage{enumerate}
\usepackage{gensymb}
\usepackage{subfigure}
\usepackage{multirow}
\usepackage{booktabs}
\usepackage{subfigure}
\graphicspath{{Images/}}

\newtheorem{proposition}{Proposition}

\newcommand{\R}{{\mathbb R}}

\newcommand{\mb}{\mathbb}
\newcommand{\mc}{\mathcal}
\newcommand{\mt}{\mathrm}
\newcommand{\bs}{\boldsymbol}

\newcommand{\scal}{0.8}

\newcommand{\infimum}{\mathop{\mathrm{inf}}}

\title[Data-Driven Robust Predictive Control for Energy Efficient Buildings and Districts]{The Power of Diversity: Data-Driven Robust Predictive Control for Energy Efficient Buildings and Districts}

\author{Georgios Darivianakis, Angelos Georghiou, Roy S. Smith and John Lygeros
	}%
	\thanks{This research was partially funded by CTI within the SCCER FEEB\&D, the Swiss National Science Foundation under the project IMES and the European Commission under the project Local4Global.} 
	\thanks{The authors are with the Automatic Control Laboratory, Department of Electrical Engineering and Information Technology, ETH Zurich, 8092 Zurich, Switzerland.
		{\tt\small \{gdarivia, angelosg, rsmith, lygeros\}@control.ee.ethz.ch}}
	
\begin{document} 
\maketitle

\begin{abstract}
	The cooperative energy management of aggregated buildings has recently received a great deal of interest due to substantial potential energy savings. These gains are mainly obtained in two ways: $ (i) $ Exploiting the load shifting capabilities of the cooperative buildings; $ (ii) $ Utilizing the expensive but energy efficient equipment that is commonly shared by the building community (e.g., heat pumps, batteries and photovoltaics). Several deterministic and stochastic control schemes that strive to realize these savings, have been proposed in the literature. A common difficulty with all these methods is integrating knowledge about the disturbances affecting the system. In this context, the underlying disturbance distributions are often poorly characterized based on historical data. In this paper, we address this issue by exploiting the historical data to construct families of distributions which contain these underlying distributions with high confidence. We then employ tools from data-driven robust optimization to formulate a multistage stochastic optimization problem which can be approximated by a finite-dimensional linear program. The proposed method is suitable for tackling large scale systems since its complexity grows polynomially with respect to the system variables. We demonstrate its efficacy in a numerical study, in which it is shown to outperform, in terms of energy cost savings and constraint violations, established solution techniques from the literature. We conclude this study by showing the significant energy gains that are obtained by cooperatively managing a collection of buildings with heterogeneous characteristics.
\end{abstract}


\section{Introduction}

Approximately 20-40\% of the total energy consumption in the developed countries is attributed to the building sector, an amount that often exceeds even the industrial and transportation sectors \cite{perez2008}. Concerns about the growing environmental impact of building energy consumption, led the EU and the US government to set the target of a net zero-energy for 50\% of their commercial buildings by 2040 \cite{Sartori2012}. In this context, active building energy management has attracted considerable attention with substantial efforts to be devoted to developing sophisticated control schemes that are capable of reducing the buildings energy impact while ensuring comfortable conditions for the building users \cite{Deori2014,Aswani2012,Maasoumy2014,OlStMo13}. Nevertheless, the opportunities for large savings within individual buildings can be limited, and depend on the specific building actuation systems and construction characteristics \cite{oldewurtel2012use,SiOlCiPr11,strurzen2016}.

Further savings can be envisaged by cooperatively managing the aggregated energy demands of a collection of buildings in a district via an energy hub. The energy hub is a conceptual entity that provides the interface between the building community and the power grid by utilizing shared energy generation, conversion and storage equipment (e.g., heat pumps, batteries and photovoltaics) \cite{geidl2007energy}. In this setting, significant energy gains can be obtained by exploiting the \emph{diversity} of the available energy sources, equipment and building characteristics. However, the main body of the literature separates the optimal control of the energy hub from that of the buildings. A number of papers treat the building energy demands as exogenous signals, which are typically estimated using building simulation environments such as EnergyPlus \cite{EnergyPlus}. These studies focus on the control of the devices within the energy hub, employing either deterministic \cite{evins2014new,fabrizio2009hourly} or stochastic \cite{parisio2012robust} formulations. 

A large body of literature addresses the building control problems using deterministic schemes \cite{OlStMo13,oldewurtel2012use,SiOlCiPr11,strurzen2016}, as they are scalable and suitable for problems with long prediction horizons. However, deterministic schemes suffer from frequent constraints violations due to their inability to handle the system disturbances \cite{Bertsimas2011}. Stochastic schemes can potentially address this issue but very often they are either not practically scalable (e.g., see the discussion in \cite{zhang2014sample} on the application of the scenario approach \cite{dupavcova2003scenario,Calafiore2006} to building control problems), or their reliability strongly depends on the exact knowledge of the distributional characteristics of the system disturbances \cite{Shapiro2011}. Poor estimation of these underlying distributions can lead to significant performance deterioration \cite{Nemirovski2006}.

To tackle this issue, distributionally robust methods have recently been introduced in the literature \cite{Wiesemann2014a}. Instead of considering an exact distribution, these methods account for all distributions in a family which shares only a few structural parameters, such as moments and/or support information \cite{VaParys13,LaGhBhJo03}. In this context, the historical data are only partially exploited to obtain a rough estimate of these structural parameters. On the other hand, purely data-driven approaches which systematically exploit the historical data information to determine appropriate families of distributions, have also been suggested \cite{BeGuKa13,Delage2010}. These methods provide the tools to reformulate an originally infinite dimensional problem to a finite dimensional convex semi-definite optimization problem (SDP). Typically, the number of constraints in this problem depends on the historical data size, which limits the scalability of the method when dealing with a large data set. Although polynomial algorithms exist for solving SDPs \cite{Wolkowicz2012}, these algorithms are computationally demanding limiting the applicability of these methods to small problem instances. 

Our goal is to develop a data-driven stochastic control scheme that is capable of cooperatively operating the energy hub and the district buildings. This paper extends the preliminary work in \cite{DaGeSmLy15}, providing a refined robust approach on handling the system disturbances by systematically exploiting the available historical realizations of these stochastic processes. In particular:
\begin{enumerate}
	\item We propose a data-driven approach that exploits the historical data to train linear models of the exogenous disturbances, and construct families of distributions that encompass the true disturbance distributions with high confidence.
	
	\item We exploit the structure of these families of distributions to formulate a robust multistage stochastic optimization problem that minimizes the wost-case expected energy costs of the system. The size of the resulting optimization problem is independent of the historical data size. We approximate this infinite dimensional problem by a finite-dimensional linear program that scales polynomially with respect to the prediction horizon length, and more importantly, unlike the SDP approaches cited above, it can effectively be solved for large scale systems. 
	
	\item We demonstrate in an extensive numerical study the efficacy of the proposed method which is shown to outperform the optimally tuned deterministic equivalent in terms of energy consumption and constraint violations. We observe that higher cost benefits are obtained by merging buildings with dissimilar operation plans, rather than diverse construction characteristics.
\end{enumerate}

The paper is organized as follows. In Section~\ref{sec::modelling}, we review in a more compact way, the modeling approach presented in \cite{DaGeSmLy15}. The main contributions of this paper are summarized in Sections~\ref{sec::distModel} and \ref{sec::probForm} where the developed data-driven distributionally robust methods, and the techniques associated with the derivation of a tractable approximation to the infinite dimensional stochastic optimization problem, are discussed. We conclude this paper with an extensive numerical study performed in Section \ref{sec::Numerical_results}. The proofs of the propositions can be found in the Appendix.

\textbf{Notation:}
All random vectors appearing in this paper are defined on an abstract probability space $(\Omega, \mc F, \mb P)$, where $\mb E(\cdot)$ denotes the expectation operator with respect to $\mb P$. Random vectors are represented in boldface, while their realizations are denoted by the corresponding symbols in normal font. For given matrices $ (A_1,\ldots,A_m) $, we define $ A := \text{diag}(A_1,\ldots,A_m) $ as the block-diagonal matrix with elements $ (A_1,\ldots,A_m) $ on its diagonal. Given vectors $ (v_1,\ldots,v_m) $, $ v_i \in \mb R^{k_i} $, we define $ [v_1,\ldots,v_m] := [v_1^\top,\ldots,v_m^\top]^\top \in \mb R^{k} $ with $ k = \sum_{i=1}^{m}k_i $, as their vector concatenation. We denote by $\bs{1}$ and $ \bs 0 $ the vectors with components all one and zero, respectively. The dimension of the corresponding vectors and the vector concatenations, will be clear from the context.

\section{System modelling}\label{sec::modelling}
In this section, we describe the energy hub and building dynamics using discrete time, bilinear models affected by stochastic exogenous disturbances. We assume that these disturbances evolve according to stochastic processes $\{\bs{\xi}_t\}_{t\in\mc T}$, where $\mc T = \{1,\ldots,T\}$, and $T$ is the length of the horizon considered. The vector $\bs{\xi}_t$ encompasses all the stochastic processes affecting the energy hub and building dynamics, such as the ambient temperature, solar radiation and internal gains of the buildings.  

\subsection{Energy hub dynamics}
\begin{figure}[t]
	\centering
	\psfrag{x1}[c][r][\scal][0]{Electricity, $ \bs d_{t,\text{elec}} $}
	\psfrag{x2}[c][r][\scal][0]{Cooling, $ \bs d_{t,\text{cool}} $}
	\psfrag{x3}[c][r][\scal][0]{Heating, $ \bs d_{t,\text{heat}} $}
	\psfrag{x4}[c][c][\scal][0]{Battery}
	\psfrag{x6}[c][c][\scal][0]{Chiller}
	\psfrag{x7}[l][l][\scal][0]{Photovoltaics}
	\psfrag{x8}[c][c][\scal][0]{Heat pump}
	\psfrag{x9}[c][c][\scal][0]{Boiler}
	\psfrag{y1}[c][l][\scal][0]{Electrical}
	\psfrag{y11}[c][c][\scal][0]{grid}
	\psfrag{y2}[c][c][\scal][0]{Gas grid}
	\psfrag{y3}[c][c][\scal][0]{Energy Hub}
	\psfrag{y4}[c][c][\scal][0]{District}
	\psfrag{p1}[c][c][\scal][0]{$ \bs p_{t,\text{el}}^{\text{in}} $}
	\psfrag{p2}[c][c][\scal][0]{$ \bs p_{t,\text{el}}^{\text{out}} $}
	\psfrag{p3}[c][c][\scal][0]{$ \bs p_{t,\text{gas}}^{\text{in}} $}
	\includegraphics[width = 0.6\textwidth]{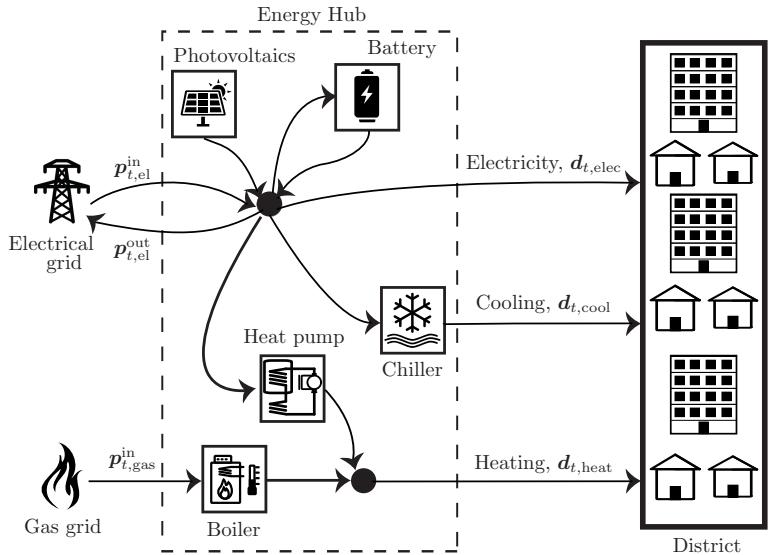}
	\caption{Heating, cooling and electricity network of a district.}
	\label{fig::enHubEx}
\end{figure}

We visualize the energy hub as a $ \emph{conceptual} $ entity that houses and interconnects a number of conversion, storage and production devices that are shared by the building community. An energy hub essentially provides the interface between the energy grid and the building community. As depicted in the illustrative example of Fig.~\ref{fig::enHubEx}, the energy hub is capable of $ (i) $ purchasing electricity and gas to meet the electricity, cooling and heating demand of the building community, and $ (ii) $ selling the electricity produced by the photovoltaics or stored in the battery to maximize profit. 
In the following, we define the set $ \mc K $ to include the energy hub devices, and the sets $\mc G$ and $\mc O$ to include the energy streams from the grid and building community, respectively. In the example of Fig.~\ref{fig::enHubEx}, $ \mc K = \{$Battery, Chiller, Heat pump, Boiler, Photovoltaics$\} $, $ \mc G  = \{$Electrical grid, Gas grid$\} $ and $ \mc O = \{$Electricity, Cooling, Heating$\} $. 

We model every device $ i \in \mc K $ using linear dynamics and constraints, as follows:
\begin{equation}\label{hub1}
	\begin{array}{l}
		\displaystyle\bs{x}_{t+1,i}  = A_i \, \bs{x}_{t,i} + B_i \,\bs{u}_{t,i} + C_i \, \bs{\xi}_{t},\\[1ex]
		F_{x,i}\, \bs x_{t,i} + F_{u,i}\, \bs u_{t,i} + F_{\xi,i} \, \bs \xi_{t}\le h_i .
	\end{array}
\end{equation}
The vectors $ \bs x_{t,i} $ and $ \bs u_{t,i} $ denote the hub device internal states and control inputs, respectively. The matrices $ A_i $, $ B_i $ and $ C_i $ have appropriate dimensions and can be derived from the device characteristics. Finally, we assume linear operational constraints captured by the matrices $ F_{x,i} $, $ F_{u,i} $, $ F_{\xi,i} $ and $ h_i $. As discussed in \cite{geidl2007energy,fabrizio2009hourly}, linear approximations for the energy hub devices are reasonable at this level of abstraction. Examples of such models for the numerical study of \mbox{Section \ref{sec::Numerical_results}} are given in Appendix \ref{app::enHub}. 

The energy balancing nodes are used to model the interconnection of the energy hub devices, as follows:
\begin{equation}\label{hub2}
	H_{p}\, \bs p_t + H_{u}\, \bs u_t  + H_{d}\,\bs d_t = \bs 0,
\end{equation}
where the vector $ \bs p_t = [\bs p_{t,1}^\text{in},\bs p_{t,1}^\text{out},\ldots,\bs p_{t,|\mc G|}^{\text{in}},\bs p_{t,|\mc G|}^\text{out}] $ contains the non negative decision variables of the grid energy streams, with $ \{\bs p_{t,i}^\text{in}\}_{i\in \mc G} $ and $ \{\bs p_{t,i}^\text{out}\}_{i\in \mc G} $ denoting the energy purchased from, and sold to the grid, respectively. In a similar way, the vector $ \bs d_t = [\bs d_{t,1},\ldots,\bs d_{t,|\mc O|}] $ concatenates the demands of the building community, and $ \bs u_t = [\bs u_{t,1},\ldots,\bs u_{t,|\mc K|}] $ captures the decision variables of the energy hub devices. The matrices $ H_p $, $ H_u $ and $ H_d $ have proper dimensions and model the power flows affecting the respective balancing node.


\subsection{Building dynamics}
We model the building dynamics using bilinear state space models, motivated by the resistance-capacitance models of \cite{sturzenegger2014brcm}. The accuracy of these bilinear models was validated against established building simulation software \cite{EnergyPlus}, and real buildings \cite{strurzen2016}. We denote by $ \mc B $ the set of district buildings and for each building $ i \in \mc B $ we assume a bilinear model that captures the temperature evolution of the rooms, walls, ceiling and floors, as follows,
\begin{equation} \label{building0}
	\displaystyle\bs{x}_{t+1,i}  = A_i \bs{x}_{t,i} + \big(B_i + \bs{x}_{t,i}^\top E_{i}\big) \bs{u}_{t,i} + \big( D_i + \bs v_{t,i}^\top C_i \big) \bs{\xi}_{t},
\end{equation}
where $\bs u_{t,i}$ contains the inputs to these building actuation systems which are not coupled with the disturbances (e.g., radiators, thermally activated building structures (TABS), air handling unit (AHU), floor heating). The rest of the control inputs (e.g. , position of the blinds) are captured by $ \bs v_{t,i} $. The system matrices $ A_i $, $ B_i $, $ C_i $, $ D_{i} $ and $ E_{i} $ depend on the specific building characteristics (e.g., number of rooms, construction material, window fraction area and actuation units). 

\begin{subequations}\label{building1}
	The state to input bilinear terms in \eqref{building0} severely complicate the design of computationally tractable controllers. To address this issue, we replace the bilinear terms  $\bs{x}_{t,i}^\top E_i \bs{u}_{t,i}$, with the linear terms  $\widehat{x}_{i}^\top E_i \bs{u}_{t,i,j}$, where $\widehat{x}_{i}$ is the initial condition for the states of building $i$. Therefore, the dynamics of the approximated system are given as:
	\begin{equation}
		\begin{array}{r}
			\bs{x}_{t+1,i} = A_i \, \bs{x}_{t,i} + B_i({\widehat x_i})\,\bs{u}_{t,i} + C_i(\bs{v}_{t,i})\bs \xi_t + D_i\bs{\xi}_t,
		\end{array}
	\end{equation}
	The matrix $ B_i(\widehat{x}_i) $ and the matrix function $ C_i(\bs v_{t,i}) $, are readily constructed by the primitive description of the dynamics in \eqref{building0}. More sophisticated approaches could also be envisioned (e.g. linearizing about a state trajectory), but we do not pursue this further here for the shake of simplicity. 
	
	We consider linear operational constraints (e.g., radiator and AHU limitations), as follows:
	\begin{equation}\label{building00}
		F_{x,i}\, \bs x_{t,i} + F_{u,i}\, \bs u_{t,i} + F_{v,i}\, \bs v_{t,i} + F_{\xi,i} \, \bs \xi_{t}\le h_i\,,
	\end{equation}
	where the matrices $ F_{x,i} $, $ F_{u,i} $, $ F_{v,i} $, $ F_{\xi,i} $ and $ f_i $ are derived using the BRCM Toolbox \cite{sturzenegger2014brcm}, and are of appropriate dimensions. In addition to the operational constraints in \eqref{building00}, we consider user specified comfort ranges given as,
	\begin{equation}\label{building3}
		\text{lb}_{t,i}\leq \bs{x}_{t,i}\leq \text{ub}_{t,i},
	\end{equation}
\end{subequations} 
where $\text{lb}_{t,i}$ and $\text{ub}_{t,i}$ are lower and upper bounds. One typically considers bounds only on the room temperatures, so many of the upper/lower bounds can be assumed to be plus/minus infinity, effectively eliminating the corresponding constraints. These bounds may vary during the day to reflect occupancy patterns (e.g., office buildings temperature bounds are often relaxed during the night time since the building is empty). We will refer to constraint set \eqref{building3} as the \textit{comfort constraints} \cite{SIA2006}. 

\subsection{Coupling of buildings to the energy hub}

We model the coupling between the energy hub and the buildings with the following set of equality constraints:
\begin{equation}\label{couplingConstraints}
	\bs d_{t,j} = \sum\limits_{i \in \mc B}  \bs{\eta}^\top_{i,j}\, \bs{u}_{t,i},\quad \forall j \in \mc O,
\end{equation}
where $ \bs \eta_{i,j} $ is a $ (0,1) $-matrix that models whether the building actuation system $ \bs{u}_{t,i} $ is connected to the $ j $-th output energy source of the hub, $ \bs d_{t,j} $. 


To simplify notation, we compactly rewrite Eqns \eqref{hub1}, \eqref{hub2}, \eqref{building1} and \eqref{couplingConstraints}, as follows:
\begin{equation}\label{eq::system}
	\left.\begin{array}{@{}l}
		\bs{x} = \bs B({\widehat x})\bs{u} + \bs C(\bs v)\bs{\xi} + \bs D\bs{\xi},\\[1ex]
		\bs F_\pi\bs{\pi}  + \bs F_\xi \bs{\xi}\leq \bs{h},
	\end{array}\right.
\end{equation} 
where the concatenated vector $ \bs \pi $ is defined as $ \bs \pi = [\bs p, \bs d, \bs x, \bs u, \bs v] $. Note that all the vectors in \eqref{eq::system} are considered over the time horizon $ \mc T $, e.g., $ \bs p = [\bs p_1,\ldots,\bs p_{|\mc T|}] $, while the vectors $ \bs x $, $ \bs u $ are defined such they include the state and input variables, respectively, for both the energy hub devices and the buildings.

\section{Modelling the uncertainty}\label{sec::distModel}
To model the evolution of the stochastic process $ \bs \xi_t $ the exact characterization of its distribution is needed. However, this information is typically unavailable. Nevertheless, historical data such as past realizations of $\bs{\xi}_t$ are usually available. A simple approach is to construct an empirical distribution  using the historical data set, and then use it as a proxy for the true distribution. However, if the data set is small, then there might be several distributions that can describe these data points. As indicated in \cite{Nemirovski2006}, if one arbitrarily chooses an element from this family of distributions, e.g., the empirical distribution, then the solution of the resulting stochastic optimization problem can  differ significantly from the solution in which the true distribution is used. In the following, we adopt a robust perspective to this problem by training linear models of the exogenous disturbances and constructing families of distributions that describe the disturbance realizations during past years. We then formulate an optimization problem whose solution addresses all distributions in the constructed family.

\subsection{Dynamics of the stochastic process}\label{model}

We denote by $ \mc D $ the set that contains all the sources of disturbances appearing in our problem, e.g., solar radiation, ambient temperature, internal gains. For each $ i\in \mc D $, we model the disturbance $ \bs{\xi}_{t,i} $ as a deterministic forecast $ f_{t,i} $, plus the stochastic error term $\bs{e}_{t,i}$. To simplify notation, we omit index $ i $ in the subsequent discussion, i.e., 
\begin{subequations}\label{eq::Model}
	\begin{equation}\label{eq::distModel}
		\bs{\xi}_{t} = f_{t} + \bs{e}_{t}.
	\end{equation}
	Weather forecasts are easily accessible from national weather services (e.g., COSMO-7 of MeteoSwiss~\cite{MeteoSwiss}) while forecasts of the anticipated internal gains can be obtained from standard weekly profiles of typical building configurations (e.g., office and residential buildings~\cite{SIA2006}). The dynamical evolution of $ \bs e_{t} $, is captured by a first order autoregressive system,
	\begin{equation}\label{eq::errorModel}
		\bs{e}_{t+1} = \alpha_{t} \bs{e}_{t} +  \bs{w}_{t},
	\end{equation} 
\end{subequations}
where $\alpha_{t}$ is time-varying constant, and $\bs w_t \in W_t$ is the stochastic process governing the stochastic evolution of the error. The autoregressive evolution of $ \bs{e}_{t} $ is motivated by systematic errors between forecast and actual realization. Indeed, the forecast provided by MeteoSwiss can differ from the true realization of the uncertain parameters due to, among other things, the spatial difference between the local weather station and the building leading to error correlation over time. Additionally, $\bs{w}_{t}$ models the noise from imperfect forecasting and possible measurement noise from the sensing devices. 

Equation \eqref{eq::Model} suggests that instead of constructing the probability distribution associated with $\bs{\xi}_{t} $, we can equivalently construct the probability distribution associated with $\bs{w}_{t}$ which we denote by $ \mb P_{t}$. 
We make the following structural assumptions for the random variables $\bs{w}_{t}$: $(i)$ they are normally distributed, $(ii)$ they are mutually independent for all time stages $t\in\mathcal{T}$ and disturbances $i\in\mathcal{D}$, and $ (iii) $ they are stationary with respect to different days (i.e., $ \alpha_t $ and the distribution of $ w_t $ for the same $ t $ but for different days is the same). Essentially these assumptions imply that first order models adequately capture the correlation over time, and one can ignore the correlation between different days and disturbances. We validate these hypotheses in Section \ref{sec::Numerical_results::DistModelVerif} based on real disturbance data. 

\subsection{Data-driven uncertainty sets}\label{datadriven}

For each disturbance $i\in\mathcal{D}$, the historical data set consists of realization and forecast pairs of the form $ \{(\xi_{t}^k,\,f_{t}^k)\} $ with $ t $ denoting the time of the $ k $-th day that the data was recorded. We consider $ N $ records of forecast and realization compatible with the sampling time used for the models of Section \ref{sec::modelling} (typically, daily records for several years sampled hourly). 

We calculate the constants $\alpha_{t}$ describing the dynamics in \eqref{eq::Model} using least squares fitting. In particular, for each $t \in \{1,\ldots,24\}$, we solve,
\begin{equation}\label{opt:alpha}
	\min_{\alpha_{t}\in\R}\;\;\sum_{k=1}^{N} \left((\xi^k_{t+1}-f^k_{t+1}) -  \alpha_{t} (\xi^k_{t}-f^k_{t})\right)^2,
\end{equation}
and using the optimal solution construct the residual data points as follows:
\begin{equation*}
	w^k_{t}:=(\xi^k_{t+1}-f^k_{t+1})-\alpha_{t} (\xi^k_{t}-f^k_{t}),~ k = 1,\ldots,N.
\end{equation*}
We denote by $\mc S_{t}:=\{w^k_{t}\}_{k=1}^{N}$ the set of residuals derived from the historical data, and subsequently denote by $\widehat{\mu}_{t}$ and $\widehat{\sigma}_{t}^2$ their empirical mean and variance, respectively.

We construct a family of distributions, $\mc P_{t}$, that are compatible with the residuals of our historical data $\mc S_{t}$. The set $\mc P_{t}$ is sometimes referred as the ambiguity set in the robust optimization literature \cite{Wiesemann2014a,BeGuKa13}. We consider an ambiguity set, $ \mc P_{t} $, of the following form:
\begin{subequations}\label{eq::ambSetBounds}
	\begin{equation}
		\begin{array}{r@{\,}l}
			\mathcal{P}_{t} = \{ \mc{N}( \mu_{t},\sigma_{t}^2):& \underline{ \mu}_{t} \leq  \mu_{t} \leq \overline{ \mu}_{t}, \;\underline{ \sigma}_{t}^2 \leq  \sigma_{t}^2 \leq \overline{ \sigma}_{t}^2\}.
		\end{array}
	\end{equation}
	The constants $(\underline{ \mu}_{t},\,\overline{ \mu}_{t},\,\underline{ \sigma}_{t}^2,\,\overline{ \sigma}_{t}^2)$, are selected to ensure that given $\mc S_{t}$, the true distribution, $ \mb P_{t}$ of $\bs{w}_{t}$, is an element of $\mathcal{P}_{t}$ with high probability. This is achieved by utilizing concepts from statistical hypothesis theory. In particular, the statistic associated with the chi-square hypothesis test, \mbox{$ h^{\chi}:=(N-1)\widehat{\sigma}_{t}^2/\sigma_{t}^2$}, follows the chi-square distribution with $(N-1)$ degrees of freedom \cite{LeRo06}. To ensure that  $\underline{ \sigma}_{t}^2 \leq  \sigma_{t}^2 \leq \overline{ \sigma}_{t}^2$, with probability at least  $1-\delta^{\chi}_{t}$, we set, 
	\begin{equation}\label{eq::uncert::varianceBounds}
		\begin{array}{l}
			\underline{ \sigma}_{t}^2 = {(N-1)\,\widehat{\sigma}_{t}^2}/{q^{\chi}_{(N-1)}({\delta^{\chi}_{t}/2})},\\[1ex]
			\overline{ \sigma}_{t}^2  = {(N-1)\,\widehat{\sigma}_{t}^2}/{q^{\chi}_{(N-1)} ({1-\delta^{\chi}_{t}/2})}. 
		\end{array}
	\end{equation}
	Here, $q^{\chi}_{(N-1)}(\cdot)$ denotes the quantile function of the chi-square distribution  with $(N-1)$ degrees of freedom.

	We construct the bounds of $\mu_{t}$ in a similar way. Given that the $\bs{w}_{t}$ follows a normal distribution and its variance $\sigma^2_{t}$ follows the chi-square distribution, \mbox{$h^{st}:=(\widehat{\mu}_{t} - \mu_{t})/\sqrt{\widehat{\sigma}^2_{t}/N}$}, follows a student $ t $-distribution with $N-1$ degrees of freedom, \cite{LeRo06}. To ensure that  \mbox{$\underline{ \mu}_{t} \leq  \mu_{t} \leq \overline{ \mu}_{t}$}, with probability at least  $1-\delta^{\text{st}}_{t}$, we set, 
	\begin{equation}\label{eq::uncert::meanBounds}
		\begin{array}{l}
			\overline{\mu}_{t} =  \widehat{\mu}_{t} + q_{N}^{st}({\delta^{\text{st}}_{t}/2})\sqrt{\widehat{\sigma}^2_{t}/N},\\[1ex]
			\underline{\mu}_{t} =  \widehat{\mu}_{t} - q_{N}^{st}({\delta^{\text{st}}_{t}/2})\sqrt{\widehat{\sigma}^2_{t}/N},
		\end{array}
	\end{equation}
\end{subequations}
where, $q_{N}^{st}(\cdot)$ denotes the quantile function of the Student's $ t $-distribution  with $(N-1)$ degrees of freedom.

The bounds given in \eqref{eq::ambSetBounds} explicitly determine the family of distributions, $ \mc P_{t} $. The following proposition provides the confidence, $ 1-\delta_t $, by which the true disturbance distribution, $ \mb P_t $, is a member of $ \mc P_t $. 
\begin{proposition}\label{prop::robGaus}
	Given the sample data, $\mc S_{t}$, let $\normalfont \bs{\mt{P}}_{\mc S_{t}} $ be the $|N|$-fold product distribution of $ \mb P_t $, then 
	\begin{equation*}
		\bs{\mt{P}}_{\mc S_{t}}\big(\mathbb{P}_{t} \in \mathcal{P}_{t}\big) \ge 1-\delta_{t}, 
	\end{equation*}
	where  $ \delta_{t} = \delta^{\text{st}}_{t}  + \delta^{\chi}_{t} $.
\end{proposition}

Let us now define the multivariate stochastic process $ \bs w = [\bs w_{0,1},\ldots,\bs w_{|\mc T|,|\mc D|}] \in W $, with $ W = W_{0,1}\times \ldots \times W_{|\mc T|,|\mc D|} $. We characterize the joint distribution $ \mb P $ of $ \bs w $, as follows,
\begin{equation*}
	\begin{array}{r@{\,}l}
		\mc P = \Big\{ \mc N (\bs \mu, \bs \Sigma):  & \textrm{ with }  \bs \mu = ( \mu_{1,1},\ldots, \mu_{|\mc T|,|\mc D|})\\
		& \hspace*{-5mm}\textrm{ and } \bs \Sigma = \mt{diag}( \sigma_{1,1}^2,\ldots, \sigma_{|\mc T|,|\mc D|}^2)\\
		&\hspace*{-13mm} \textrm{ s.t.} \mb P_{1,1}  \in \mc P_{1,1}, \ldots, \mb P_{|\mc T|,|\mc D|} \in \mc P_{|\mc T|,|\mc D|} \Big\}.
	\end{array}
\end{equation*}
We close this section, by compactly rewriting the disturbance modelling equations over the horizon $ T $, as follows: 
\begin{equation*}
	\bs \xi = \bs H(\widehat{e}) \bs w,
\end{equation*}
where the matrix $ \bs H(\widehat{e}) $ is readily constructed from \eqref{eq::Model} and depends on the forecast, and the vector $ \widehat{e} = [e_{0,1},\ldots,e_{0,|\mc D|}] $ containing the current error disturbance measurements.

\section{Problem formulation} \label{sec::probForm}
{Our objective is to minimize the worst-case expected cost of energy purchased from the grid by the building community, over a finite horizon.} This must be achieved while satisfying the dynamics and constraints of the devices and the buildings in the system. This problem can be  formulated as a multistage stochastic linear program, as follows: 
\begin{equation}\label{problem_bilinear}
	\begin{array}{l@{\,}l}
		\textrm{min}\;\;& \sup\limits_{\mb P \in \mc P}\displaystyle \;\mb E_{\mb P} \left( \bs c^\top \bs{p}\right)  \\[2ex]
		\textrm{ s.t.}& \bs p, \bs{u}\in \mathcal{C},\;v \in \mathcal{R},\; \bs \pi = [\bs p,\bs d,\bs x,\bs u,v],\\[1ex]
		& \bs \xi = \bs H(\widehat{e}) \bs w,\\[1ex]
		&\bs{x} = \bs B({\widehat x})\bs{u} + \bs C(v)\bs{\xi} + \bs D\bs{\xi},\\[1ex]
		&\inf\limits_{\mb P \in \mc P} \mb P \Big( \bs F_\pi\bs{\pi}  + \bs F_\xi \bs{\xi}\leq \bs{h}\Big) \geq 1-\epsilon,
	\end{array}
\end{equation}
where $ \bs c_t = [ c_{t,1}^\text{in}, c_{t,1}^\text{out},\ldots, c_{t,|\mc G|}^{\text{in}}, c_{t,|\mc G|}^\text{out}]  $ contains the time varying prices of the grid energy streams. Note that the prices $ \{c_{t,i}^\text{in}\}_{i\in \mc G} $ for purchasing energy from the grid are positive scalars, while $ \{c_{t,i}^\text{out}\}_{i\in \mc G} $ are negative. Furthermore, we select a chance constraints formulation with worst-case violation probability $ \epsilon $. The decision variables $ v $ that are coupled with the disturbances $ \bs \xi $, are allowed to take values in the generic finite-dimensional vector space $ \mc R $. We choose the decision variables $ \bs p, \bs u $ to be \emph{strictly causal} disturbance feedback policies. For instance, $ \bs p_t $ is the strictly causal, vector valued function of the energy purchased from the grid at time $ t\in \mc T $, defined as $ \bs p_t: W_1\times \ldots \times W_{t-1} \rightarrow \mc R $. In the following, we denote by $ \mc C $, the infinite-dimensional function space of \emph{strictly causal} disturbance feedback policies. Finally, the expected cost of the objective has been formulated to account for the worst-case multivariate distribution $ \mb P $ in the ambiguity set $ \mc P $.

\subsection{Constraint relaxation}

The proposed control architecture is based on a receding horizon implementation of the system. In this setting, the mismatch between the building prediction model and real system, can lead to comfort bounds violations. To address infeasible instances of Problem \eqref{problem_bilinear} due to this issue, we relax the comfort constraints \eqref{building3}, as follows:
\begin{equation*}
	\left.\begin{array}{r@{\,}l}
		\max\{\text{lb}_{t,i} - \bs{x}_{t,i},0, \bs{x}_{t,i} - \text{ub}_{t,i}\} &\leq \bs{s}_{t,i},
	\end{array}\right.
\end{equation*}
where we refer $ \bs s_{t,i} $ as the slack variable. To this end, we rewrite the compactly formulated inequality \eqref{eq::system}, as follows:
\begin{equation}\label{slack}
	\left.\begin{array}{@{}l}
		\bs F_\pi\bs{\pi}  + \bs F_\xi \bs{\xi}-\bs s \leq \bs{h},
	\end{array}\right.
\end{equation}
where $ \bs s \in \mc C $ denotes the concatenated slack variable vectors over the horizon $ T $. Notice that the slack variables are only introduced to deal with building model mismatch which can be an issue on a receding horizon implementation. On the other hand, the chance constraint formulation in Problem \eqref{problem_bilinear}, is primarily used to address in a probabilistic fashion the extreme realization of the uncertain parameters $ \bs w $.

\subsection{Optimization over linear feedback policies}\label{subsec::solApprox}

The optimal solution of Problem \eqref{problem_bilinear} remains intractable due to the infinite dimensional structure of its decision variables. However, a tractable approximation of Problem \eqref{problem_bilinear} can be obtained by restricting the decision variables to the finite dimensional space of affine policies denoted as $ \mc C_{\text{aff}} $. A strictly causal affine policy, e.g. for the variables $ p_t $ is given by
\begin{equation}\label{eq::linPol}
	\bs p_t = p_{0,t} + \displaystyle \sum_{s = 1}^{t-1} P_{s,t} \,\bs w_s,
\end{equation}
where $ p_{0,t} \in \mc R $, and matrices $ P_{s,t} $ of appropriate dimensions. We refer to \eqref{eq::linPol} as the strictly causal affine decision rule (ADR) \cite{Georghiou:10}. We refer to the linear policy in \eqref{eq::linPol} which disregards the disturbance history (e.g., $ \bs p_t = p_{0,t} $), as open loop policy (OLP). The non-adaptive nature of the OLP policy provides an even more conservative {controller parametrization} with the benefit of a considerably smaller number of optimization variables.

In this context, the approximated variant of Problem \eqref{problem_bilinear} is given as follows:
\begin{equation}\label{problem_bilinearv2}
	\begin{array}{l@{\,}l}
		\textrm{min}\;\;& \sup\limits_{\mb P \in \mc P}\displaystyle \;\mb E_{\mb P} \left( \bs c^\top \bs{p}+\gamma \bs{1}^\top\bs{s}\right)  \\[2ex]
		\textrm{ s.t.}& \bs p, \bs{u}, \bs s\in \mathcal{C}_{\text{aff}},\;v \in \mathcal{R},\; \bs \pi = [\bs p,\bs d,\bs x,\bs u,v],\\[1ex]
		&\bs \xi = \bs H(\widehat{e}) \bs w,\\[1ex]
		&\bs{x} = \bs B({\widehat x})\bs{u} + \bs C(v)\bs{\xi} + \bs D\bs{\xi},\\[1ex]
		&\inf\limits_{\mb P \in \mc P} \mb P \Big( \bs F_\pi\bs{\pi}  + \bs F_\xi \bs{\xi}-\bs s\leq \bs{h}\Big) \geq 1-\epsilon,
	\end{array}
\end{equation}
where the additional term $\gamma \bs{1}^\top\bs{s}$ in the objective penalizes the constraint violations in \eqref{slack}, with parameter $\gamma\in\R_{+}$. A discussion for appropriate values of $ \gamma $ can be found in \cite{Kerrigan2000}.

\subsection{Chance constraint approximation}

The main body of the literature exploits the structure of the ambiguity sets to propose semi-definite \cite{VaParys13,Delage2010} and second-order cone \cite{calaElGha06} reformulations of the distributionally robust chance constraints. Although convex problems with conic constraints are generally tractable, they are computationally demanding for large systems. An alternative approach is to construct a set $ \widehat{W} \subseteq W $ such that the feasible region of the corresponding robust constraint is a subset of the feasible region of the  distributionally robust chance constraint, i.e., 
\begin{subequations}\label{cond::chCon}
	\begin{align}
		\text{if } \bs F_\pi\bs{\pi}  + \bs F_\xi \bs{\xi}-\bs s\leq \bs{h},\, \forall \, \bs w\in \widehat{W}, \label{cond::chCon::a} \\
		\text{then } \inf\limits_{\mb P \in \mc P} \mb P \Big( \bs F_\pi\bs{\pi}  + \bs F_\xi \bs{\xi}-\bs s\leq \bs{h}\Big) \geq 1-\epsilon.
	\end{align}
\end{subequations}
Such methods are discussed in \cite{BeGuKa13} where the authors compute $\widehat W$ described by non-linear constraints. Nevertheless, these approaches typically lead to semi-definite reformulations of the robust constraint \eqref{cond::chCon::a}. To tackle this issue, we resort to a more stringent condition, in which we require that at least, $ 1-\epsilon $, of the probability mass of each distribution, $ \mb P $, in the ambiguity set $ \mc P $, is contained in $ \widehat W $. To this end, we construct for every $ i \in \mc D $ and $ t\in \mc T $, the compact convex set $ \widehat{W}_{t,i} $, as follows:
\begin{subequations}\label{eq::setW}
	\begin{equation}
		\begin{array}{@{}r@{\,}l}
			\widehat{W}_{t,i} = \Big\{\bs w_{t,i}\;|\;&\displaystyle  \bs w_{t,i} \ge \underline{ \mu}_{t,i} - \Phi^{-1}\Big(1-\underline{\beta}_{t,i}\Big) \, \overline{\sigma}_{t,i}\,,\\[1ex]
			&\displaystyle  \bs w_{t,i} \le  \overline{ \mu}_{t,i} + \Phi^{-1}\Big(1-\overline{\beta}_{t,i}\Big) \, \overline{\sigma}_{t,i}~ \Big\}.
		\end{array}
	\end{equation}
	where the constants $ \underline{ \mu}_{t,i} $, $ \overline{ \mu}_{t,i} $, and $ \overline{\sigma}_{t,i} $ are given in \eqref{eq::uncert::meanBounds} and $ \eqref{eq::uncert::varianceBounds} $, respectively. Moreover, $ \Phi^{-1}(\cdot) $ denotes the inverse cumulative normal distribution function, and $ \underline{\beta}_{t,i} $ and $ \overline{\beta}_{t,i} $ are positive constants chosen as,
	\begin{equation}
		\displaystyle \sum_{t \in \mc T} \sum_{i \in \mc D} \Big(\underline{\beta}_{t,i} + \overline{\beta}_{t,i}\Big) = \epsilon.
	\end{equation}
	Finally, we set $ \widehat{W} $ to be,
	\begin{equation}
		\widehat{W} = \widehat{W}_{1,1} \times \ldots \times \widehat{W}_{|\mc T|,|\mc D|}.
	\end{equation}
\end{subequations}

\begin{proposition}\label{prop::strCon}
	Let $ \widehat{W} $ defined in \eqref{eq::setW}. Then, the following probabilistic guarantee holds:
	\begin{equation*}
		\infimum_{\mb{P}  \in \mc{P} } \mb{P}(\bs w \in \widehat{W}) \ge 1-\epsilon.
	\end{equation*}
\end{proposition}
Notice that satisfying \eqref{cond::chCon::a} together with the condition in Proposition \ref{prop::strCon} is more stringent than the actual chance constraint. In other words, 
\begin{equation*}
	\begin{array}{l}
		\text{if } \inf\limits_{\mb{P}  \in \mc{P} } \mb{P}(\bs w \in \widehat{W}) \ge 1-\epsilon, \text{ and \eqref{cond::chCon::a} holds, } \\
		\text{then } \inf\limits_{\mb P \in \mc P} \mb P \Big( \bs F_\pi\bs{\pi}  + \bs F_\xi \bs{\xi}-\bs s\leq \bs{h}\Big) \geq 1-\epsilon.
	\end{array}
\end{equation*}
Imposing this additional condition, we further restrict the feasible region with the benefit of gaining computational tractability due to the simple hyperrectangular structure of $ \widehat{W} $. This, in turn, allows us to address the large scale problems examined in this paper.

\subsection{Objective function reformulation}

We employ the epigraph representation to equivalently rewrite the worst-case expectation in the objective function of Problem \eqref{problem_bilinearv2}. In particular, the linearity of the objective function allows us to replace the expected value with:
\begin{equation}\label{eq::expObj}
	\sup\limits_{\mb P \in \mc P} \mb E_{\mb P} \left( \bs c^\top \bs{p}+\gamma \bs{1}^\top\bs{s}\right) = \sup\limits_{\mb P \in \mc P} \left(\bs c^\top \widetilde{\bs{p}}+\gamma \bs{1}^\top\widetilde{\bs{s}}\right) = \tau,
\end{equation}
which can equivalently be written as,
\begin{equation}\label{eq::expRep}
	\bs c^\top \widetilde{\bs{p}}+\gamma \bs{1}^\top \widetilde{\bs{s}} \le \tau, \; \forall \bs \mu \in [\underline{\bs \mu},\overline{\bs \mu}],
\end{equation}
where $ \widetilde{\bs{p}} $, and $ \widetilde{\bs{s}} $ are derived from \eqref{eq::linPol} by replacing the stochastic variable $ \bs w $ by its expected value $ \bs \mu $ (e.g., $ \widetilde{\bs{p}}_t = p_{0,t} + \sum_{s = 1}^{t-1}P_{s,t} \,\bs \mu_{t} $). Notice that the epigraph representation in \eqref{eq::expRep} is exact with the scalar variable $ \tau $ replacing the original worst-case objective.

We conclude this section by providing the final approximation of Problem \eqref{problem_bilinear} as follows: 
\begin{equation}\label{problem_bilinearv3}
	\begin{array}{l@{\,}l}
		\textrm{min}\;\;& \tau \\
		\textrm{ s.t.}& \bs p, \bs{u}, \bs s\in \mathcal{C}_{\text{aff}},\;v \in \mathcal{R},\; \bs \pi = [\bs p,\bs d,\bs x,\bs u,v],\\[1ex]
		& \bs c^\top \widetilde{\bs{p}}+\gamma \bs{1}^\top \widetilde{\bs{s}} \le \tau, \; \forall \bs \mu \in [\underline{\bs \mu},\overline{\bs \mu}],\\[1ex]
		&\left.\begin{array}{@{}l}
			\bs \xi = \bs H(\widehat{e}) \bs w,\\[1ex]
			\bs{x} = \bs B({\widehat x})\bs{u} + \bs C(v)\bs{\xi} + \bs D\bs{\xi},\\[1ex]
			\bs F_\pi\bs{\pi}  + \bs F_\xi \bs{\xi}-\bs s\leq \bs{h},
		\end{array}\right \rbrace \forall \bs w \in \widehat{W}.
	\end{array}
\end{equation}
Problem \eqref{problem_bilinearv3} still retains its infinite structure involving a continuum space of decision variables and constraints. However, as indicated in \cite{bental2004ars,Georghiou:10}, Problem \eqref{problem_bilinearv3} can be reformulated into a linear optimization problem by employing traditional robust optimization techniques that use duality to translate the semi-infinite structure into a finite number of linear constraints. Unlike other robust convex programs which employ semi-definite programming, this problem is computationally tractable and can be solved in very high dimensions.

\section{Numerical results}\label{sec::Numerical_results}

In this section, we perform numerical studies to assess the performance of the ADR and OLP control methodologies. We compare these methods to the classical and commonly used in practice \emph{certainty equivalence problem} (CEP) in which all random variables of Problem~\eqref{problem_bilinear} are replaced with their expected value. To make the comparison of these control strategies more transparent, we do not consider the possibility of selling energy to the grid. Moreover, we analyze the sensitivity of the considered control methods with respect to parameters deviations, and investigate the potential gains that can be obtained by cooperatively managing buildings with heterogeneous construction and operation characteristics.

\subsection{Problem configuration} \label{sec::Numerical_results::problemSetUp}

\begin{figure}[t]
	\centering
	\psfrag{x1}[c][c][\scal][0]{Electricity}
	\psfrag{x2}[c][c][\scal][0]{Cooling}
	\psfrag{x3}[c][c][\scal][0]{Heating}
	\psfrag{x4}[l][l][\scal][0]{Battery}
	\psfrag{x6}[c][c][\scal][0]{Chiller}
	\psfrag{x7}[l][l][\scal][0]{PV}
	\psfrag{x8}[c][c][\scal][0]{Heat pump}
	\psfrag{x9}[c][c][\scal][0]{Boiler}
	\psfrag{y11}[c][c][\scal][0]{Electricity}
	\psfrag{y12}[c][c][\scal][0]{grid}
	\psfrag{y2}[c][c][\scal][0]{Gas grid}
	\psfrag{y3}[c][c][\scal][0]{Energy Hub}
	\psfrag{y4}[c][c][\scal][0]{District}
	\includegraphics[width = 0.6\textwidth]{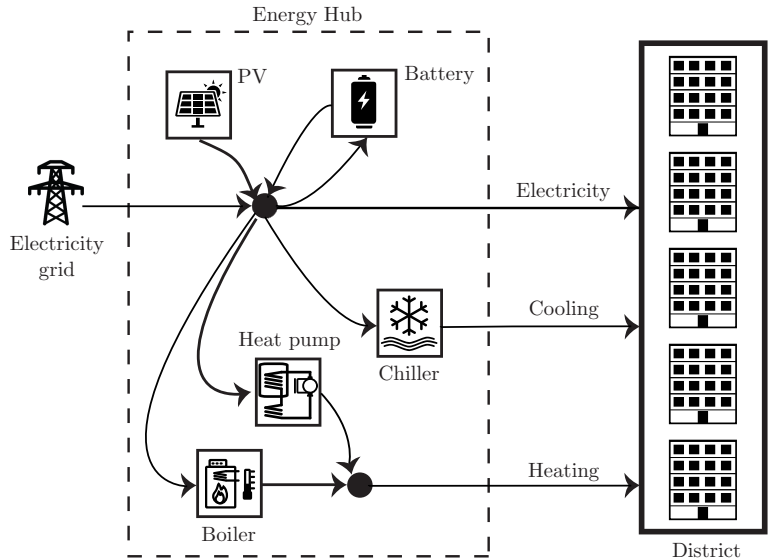}
	\caption{Energy hub configuration for numerical studies.}
	\label{fig::enHubSim}
\end{figure}

We consider districts composed of buildings with roughly the same floor area which are connected through an energy hub that comprises five devices: chiller, boiler, heat pump, photovoltaics (PV) and battery, as depicted in Fig. \ref{fig::enHubSim}. To keep the ratio between demand and supply in the energy hub, relatively constant, we linearly scale the capacities of the hub devices based on the number of buildings in our system. A detailed description of the hourly discretized dynamics and constraints governing the energy hub components can be found in the Appendix \ref{app::enHub}.

\begin{table}[b]
	\renewcommand{\arraystretch}{1.0}
	\centering
	\begin{tabular}{cccccl}
		\multicolumn{6}{c}{Building specifications}\\
		\hline
		No. & Area$(\textrm{m}^2)$  & WFA & BT & CT & \multicolumn{1}{c}{Input Devices} \\
		\hline \hline
		1 & 420   & 30\% & SP & heavy & AHU, blinds, radiator \\
		
		2 & 420   & 50\% & SP & light & AHU, blinds, TABS \\ 
		
		3 & 441   & 80\% & SA & light & AHU, blinds, TABS \\ 
		
		4 & 441   & 50\% & SA & heavy & AHU, blinds, radiator \\ 
		
		5 & 374   & 50\% & ST & heavy & AHU, blinds, radiator \\ 
		
		\hline
	\end{tabular}
	\vspace*{1mm}
	\caption{Summary of the $5$  buildings used in the simulations. }
	\label{tab:building_data}
\end{table}
The buildings have heterogeneous construction characteristics, summarized in Table \ref{tab:building_data}. In particular, we use hourly discretized building models described in  \cite{strurzen2016} characterized by the building type $ \textrm{BT}$ $\in$ $\{$\textrm{Swiss Passive (SP), Swiss Average(SA), Swiss Target (ST)}$ \} $, and construction type $\textrm{CT}$ $=$ $\{$\textrm{heavy, light}$\}$.  Each building consists of $ 5 $ rooms which are characterized by the window fraction area $\textrm{WFA}$ $=$ $\{$30\%,50\%,80\%$\}$ and their corresponding facade orientation. The dynamics and constraints of each building are generated using the BRCM toolbox \cite{sturzenegger2014brcm} to which the individual building construction details along with the specifications of the control devices (radiators, AHU, TABS, blinds) are provided.

Each of the building models in Table \ref{tab:building_data}, consists of 113 states. Model reduction techniques are used to derive a simplified, but sufficiently accurate, building model. For instance, using the \emph{bilinear balanced truncation} method described in \cite{BeDa11}, we generate a 59 state model, with a maximum absolute error of less than 0.1\degree C, during step response simulations for building 1 of Table \ref{tab:building_data}. By evaluating the Hankel matrix of this building system, we identified the time constants associated with the three most controllable and observable modes as $ 11.74 $ days, $ 4.47 $ hours and $ 4.24 $ minutes. Approximately the same time constants were identified for the other building systems in Table \ref{tab:building_data}. {The reduced order linear models are used only for prediction purposes, while the original bilinear models, given in \eqref{building0}, are employed to simulate the buildings dynamics in the closed-loop implementation}.
\begin{table}[t]
	\centering
	\renewcommand{\arraystretch}{1.0}
	\begin{tabular}{cccc}
		\hline
		Time & Cost & Winter Bounds & Summer Bounds \\
		\hline
		\hline
		05:00 - 23:00 & 0.145\text{ $ \text{CHF}/\text{kWh} $} & [21, 25]\degree C & [20, 23]\degree C \\
		23:00 - 05:00 & 0.097\text{ $ \text{CHF}/\text{kWh} $} & [15, 30]\degree C & [15, 30]\degree C\\
		\hline
	\end{tabular}
	\vspace*{1ex}
	\caption{Electricity day/night tariff variations and comfort constraints bounds. }
	\label{tab::varyingCosts}
\end{table}

We compare the performance of the ADR, OLP and CEP control designs using the metrics of purchased grid energy, room constraint violations and solution time. The cost of the energy purchased from the grid is measured in Swiss Franc (CHF), and the room constraint violations are measured in Kelvin hours (Kh). We assume time-varying electricity tariffs and comfort constraints bounds, as given in Table~\ref{tab::varyingCosts}. The same comfort bounds are used in every room and building considered in the system. The disturbance forecasts and realizations are the same as those used by the OptiControl project \cite{opticontrol2010}, for the city of Z{\"u}rich during the years 2006 and 2007. The data set from year 2006 is used to train the linear disturbance models and construct the families of distributions given in equations \eqref{opt:alpha} and \eqref{eq::ambSetBounds}. The 2007 dataset is then used to test the performance of the controllers developed based on the resulting distribution families. Finally, we choose as soft constraint penalization, $ \gamma = 10^3 $, confidence levels $ \delta^\chi_{t,i} = \delta^{\text{st}}_{t,i} = 0.01$, and constraint violation level $ \epsilon =0.01 $. We select the violation levels for the upper and lower bounds of the $ i $-th disturbance at time $ t $, as $ \overline{\beta}_{t,i} = \underline{\beta}_{t,i} ={\epsilon}/(2|\mc T||\mc D|)$.

\subsection{Disturbance model verification}\label{sec::Numerical_results::DistModelVerif}

We consider a disturbance set $ \mc D $ which comprises seven sources of uncertainty; ambient and ground temperatures, four sources of solar radiation (North, South, West and East), and building internal gains. The analysis in Section \ref{sec::distModel} requires that the residual uncertain parameters $ \bs w_t $ are normally distributed and independent over time. We verify these assumptions using the Shapiro-Wilk~\cite{RaWa11} and Pearson~\cite{LeRo06} hypothesis tests on historical, weather and occupancy, data for the year 2006.
\begin{table}[ht]
	\centering
	\renewcommand{\arraystretch}{1.0}
	\begin{tabular}{|c|cccc|}
		\hline
		Time & AT & SRS & SRE & IG \\ 
		\hline \hline
		07:00& 0.989&   0.991&  0.972 & 0.940 \\
		08:00& 0.996&   0.997&  0.996 & 0.945 \\
		09:00& 0.994&   0.995&  0.993 & 0.923\\
		10:00& 0.993&   0.965&  0.992 & 0.967\\
		11:00& 0.994&   0.990&  0.987 & 0.920\\
		12:00& 0.992&   0.941&  0.981 & 0.928\\
		13:00& 0.997&   0.996&  0.988 & 0.945\\
		14:00& 0.995&   0.995&  0.984 & 0.965\\\hline
	\end{tabular}
	\vspace*{1ex}
	\caption{Statistic of the Shapiro-Wilk normality test.}
	\label{tab::normAssump}
\end{table}

In Table \ref{tab::normAssump}, we report the values of the Shapiro-Wilk test statistic from 7:00 to 14:00 for the ambient temperature (AT), solar radiation south (SRS) and east (SRE), and internal gains (IG). We choose to present the results from 7:00 to 14:00 since these are the hours with the greatest variation in mismatch between forecasts and realizations. Values which are close to 1 indicate that the sample data are compatible with a normal distribution. In particular, the normality hypothesis is accepted at the significance level of $ 0.1 $ for the atmospheric processes, and $ 0.01 $ for the internal gains. These results are in accordance with the study in \cite{Jewson2003} where it is shown that normal distributions sufficiently capture the evolution of weather processes. By contrast, in \cite{Page2008} it is argued that Poisson distributions should be used to generate occupancy profile trajectories. Although, Poisson distributions can be approximated by Gaussians \cite{Wallace1958}, we stress that the bounds in \eqref{eq::ambSetBounds} are useful even if the underlying data are mildly non-Gaussian, as suggested in \cite[\S11.3]{LeRo06}.  

\begin{table}[ht]
	\centering
	\renewcommand{\arraystretch}{1.0}
	\begin{tabular}{|r|rrrrrr|}
		\hline Time & 07:00& 08:00& 09:00& 10:00& 11:00& 12:00 \\\hline
		07:00 & 1.00 &    0.04 &    0.02 &    0.01 &   -0.01 &    0.01 \\
		08:00 & 0.04 &    1.00 &    0.05 &    0.05 &   -0.01 &   0.00 \\
		09:00 & 0.02 &    0.05 &    1.00 &    0.06 &    0.03 &   -0.09 \\
		10:00 & 0.01 &    0.05 &    0.06 &    1.00 &    0.00 &   -0.09 \\
		11:00 & -0.01 &   -0.01 &    0.03 &    0.00 &    1.00 &   -0.03 \\
		12:00 & 0.01 &   0.00 &   -0.09 &   -0.09 &   -0.03 &    1.00 \\\hline
	\end{tabular}
	\vspace*{1ex}
	\caption{Statistic of the Pearson correlation test.}
	\label{tab::correlTest}
\end{table}
In Table \ref{tab::correlTest}, we provide the values of the Pearson test statistic from 7:00 to 12:00 for the SRS and the results are similar for the other sources of uncertainty. This correlation statistic can range between plus and minus one with values close to zero denoting probably uncorrelated data. We emphasize that for every disturbance, the hypothesis of uncorrelated data is accepted at the significance level of 0.05 which provides a strong statistical evidence for the validity of our uncertain disturbance assumptions.

\subsection{Prediction horizon selection}\label{sec::Numerical_results::sensitivityAnalysis}

We investigate the effect of the prediction horizon on the performance of the ADR, OLP and CEP control methodologies. 
\begin{figure}[t]
	\centering
	\psfrag{tL1}[c][c][\scal][0]{\textbf{Grid purchased energy}}
	\psfrag{tL2}[c][c][\scal][0]{\textbf{Room constraint violations}}
	\psfrag{yL1}[c][c][\scal][180]{Energy cost (CHF/Week)}
	\psfrag{y2}[c][c][\scal][180]{Violations (Kh/Week)}
	\psfrag{xL1}[c][c][\scal][0]{horizon (hours)}
	\psfrag{x2}[c][c][\scal][0]{horizon (hours)}
	\psfrag{tL}[c][c][\scal][0]{\textbf{Solution time}}
	\psfrag{yL}[c][c][\scal][180]{Time (sec)}
	\psfrag{xL}[c][c][\scal][0]{Horizon (hours)}
	\psfrag{MPC}[l][l][\scal][0]{\scriptsize CEP}
	\psfrag{ROB}[l][l][\scal][0]{\scriptsize OLP}
	\psfrag{LDR}[l][l][\scal][0]{\scriptsize ADR}
	\includegraphics[width = 0.6\textwidth]{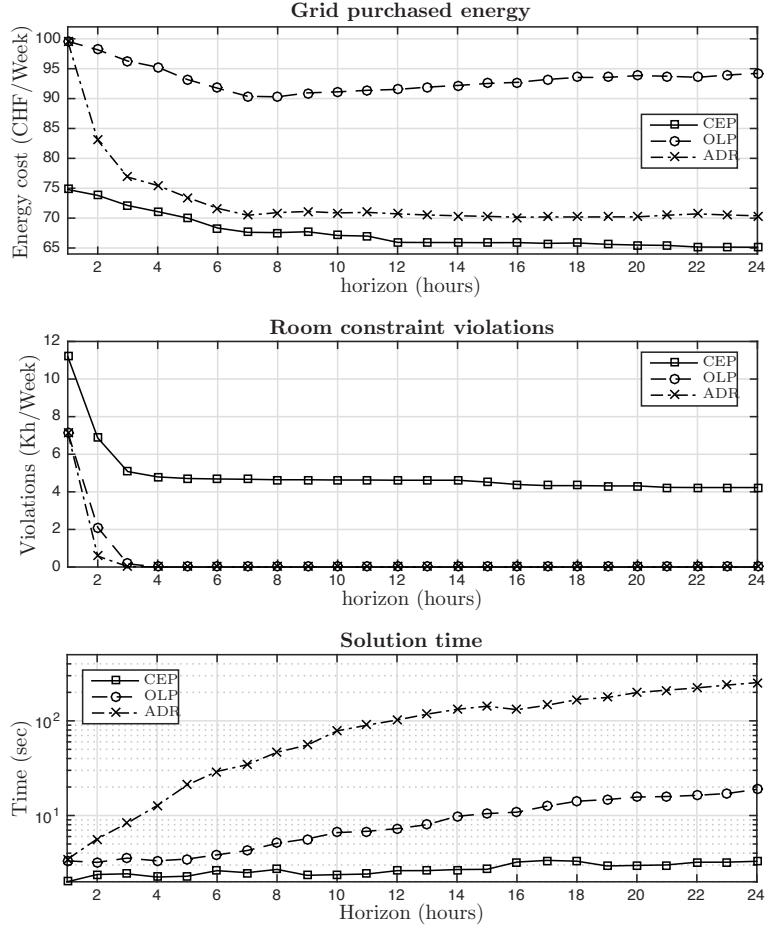}
	\caption{Purchased grid energy, room constraint violations and solution time with respect to the prediction horizon length $ T $ during a typical winter week.}
	\label{fig::PrHorSel::pareto}
\end{figure}
In Fig.~\ref{fig::PrHorSel::pareto}, we show the respective curves generated by conducting a receding horizon simulation during the first week of January $ 2007 $ for building $ 1 $ in Table~\ref{tab:building_data}. We observe that the cost of purchased grid energy associated with the CEP solution method is the least, at the expense of frequent comfort constraint violations. On the contrary, the stochastic approximations (OLP and ADR), are less cost efficient but lead to many fewer constraint violations. 

We can identify several factors associated with the selection of a suitable prediction horizon:
\begin{itemize}
	\item[$ (i) $] The dominant time constants of the buildings were roughly identified as $ 11 $ days, $ 4 $ hours and $ 4 $ minutes. Hence, a horizon of $ T \ge 7 $ hours facilitates the controller to excite the modes of the system that determine its short term evolution. We emphasize that contrary to the ADR and CEP methods, the long horizons deteriorate the performance of the OLP method due to its inability to adapt on the growing size of the disturbance uncertainty. 
	
	\item[$ (ii) $] The six hour gap between day and night comfort bounds, given in Table \ref{tab::varyingCosts}. A prediction horizon of $ T \ge 6 $, is required for the system to anticipate the comfort bounds of the next day. In this way, it can utilize the heating capacity of the buildings and the battery storage to exploit the day-night tariff structure of the electricity prices.
	
	\item[$ (iii) $] The renewable energy production peaks around midday. A horizon of $ T \ge 6 $, is sufficient for the system to anticipate the cost-free energy that will be available from the photovoltaic units. In this way, it can efficiently utilize the battery to fully exploit the potential excess of energy. 
	
	\item[$ (iv) $] The computation time associated with the ADR and OLP control methods which increases with the prediction horizon length.
\end{itemize}
The above analysis suggests that a prediction horizon of $ T = 8 $ hours provides a reasonable trade-off between foresight of the controller and computational tractability. This value is used for the rest of this section.

\subsection{Solution method selection}\label{sec::Numerical_results::controlMethodSel}

To generalize our observations regarding the comparison of the ADR, OLP and CEP control methodologies. The simulation experiment is extended to a district composed of the $ 5 $ buildings summarized in Table~\ref{tab:building_data}. The system is simulated in a receding-horizon fashion using data realizations of $12$ consecutive weeks (restarting at the beginning of each week) for the winter and summer periods of $2007$, starting January 1st and June 29th, respectively. 
\begin{table}[t]
	\centering
	\renewcommand{\arraystretch}{1.0}
	\begin{tabular}{c|ccc}
		\multicolumn{1}{c}{~}&\multicolumn{3}{c}{\emph{\textbf{Winter}}} \\
		\hline
		\multicolumn{1}{c|}{Method} & \multicolumn{1}{c}{Cost (p.u.)} & \multicolumn{1}{c}{Violations (Kh/Week)} & Basis  (CHF/Week)\\
		\hline 
		\hline 
		CEP & $(1.00,\; 0.21)$ & $(4.31,\; 0.68)$  & 327.56\\
		OLP & $(1.34,\; 0.19)$ & $(0.00,\; 0.00) $ & \\
		ADR & $(1.07,\; 0.21)$ & $(0.03,\; 0.01) $ &  \\ 
		\hline
	\end{tabular}
	
	\vspace*{5mm}
	\begin{tabular}{c|ccc}
		\multicolumn{1}{c}{~}&\multicolumn{3}{c}{\emph{\textbf{Summer}}} \\
		\hline
		\multicolumn{1}{c|}{Method} & \multicolumn{1}{c}{Cost (p.u.)} & \multicolumn{1}{c}{Violations (Kh/Week)} & Basis  (CHF/Week)\\
		\hline 
		\hline 
		CEP &  $(1.00,\; 0.17)$ & $(2.23,\; 0.14)$ & 64.19 \\
		OLP &  $(1.07,\; 0.18)$ & $(0.03,\; 0.01)$ &  \\
		ADR &  $(1.03,\; 0.16)$ & $(0.06,\; 0.01)$ & \\ 
		\hline
	\end{tabular}
	\vspace*{2mm}
	\caption{Receding horizon performance results.}
	\label{tab:results}
\end{table}

We calculated the purchased energy cost and room constraint violations for each one of these 12 weeks, and we report the results in Table~\ref{tab:results}. The table entries correspond to the \emph{(empirical mean, empirical standard deviation)} over these 12 weeks scenarios. Notice that we present the energy costs in the per unit (p.u.) system using as base values the mean costs occurring for the CEP method during the winter and summer period, respectively. Once again, we observe that the CEP solution method is the most cost effective at the expense of significant comfort constraint violations. The ADR and OLP control methodologies result in higher cost of purchased grid energy but results in fewer constraint violations. Notice that the ADR approximation is considerably more cost efficient than the OLP, while it achieves the same level of constraint violation. This is attributed to the nature of the ADR method which takes into account the potential adaptation of future decisions based on the realizations of the disturbance variables.

\begin{figure}[ht]
	\psfragscanon
	\psfrag{x1}{\Large{(\degree C)}}
	\begin{center}
		\psfrag{tL1}[c][c][\scal][0]{\textbf{Room temperature}}
		\psfrag{yL1}[c][c][\scal][180]{Temperature (\degree C)}
		\psfrag{xL1}[c][c][\scal][0]{hours}
		\psfrag{tL2}[c][c][\scal][0]{\textbf{Grid purchased energy}}
		\psfrag{yL2}[c][c][\scal][180]{Power (kW)}
		\psfrag{xL2}[c][c][\scal][0]{hours}
		\psfrag{tL3}[c][c][\scal][0]{\textbf{Grid purchased energy}}
		\psfrag{yL3}[c][c][\scal][180]{Energy cost (p.u.)}
		\psfrag{xL3}[c][c][\scal][0]{hours}
		\psfrag{MPC}[l][l][\scal][0]{\scriptsize CEP}
		\psfrag{ROB}[l][l][\scal][0]{\scriptsize OLP}
		\psfrag{LDR}[l][l][\scal][0]{\scriptsize ADR}
		\includegraphics[width = 0.6\textwidth]{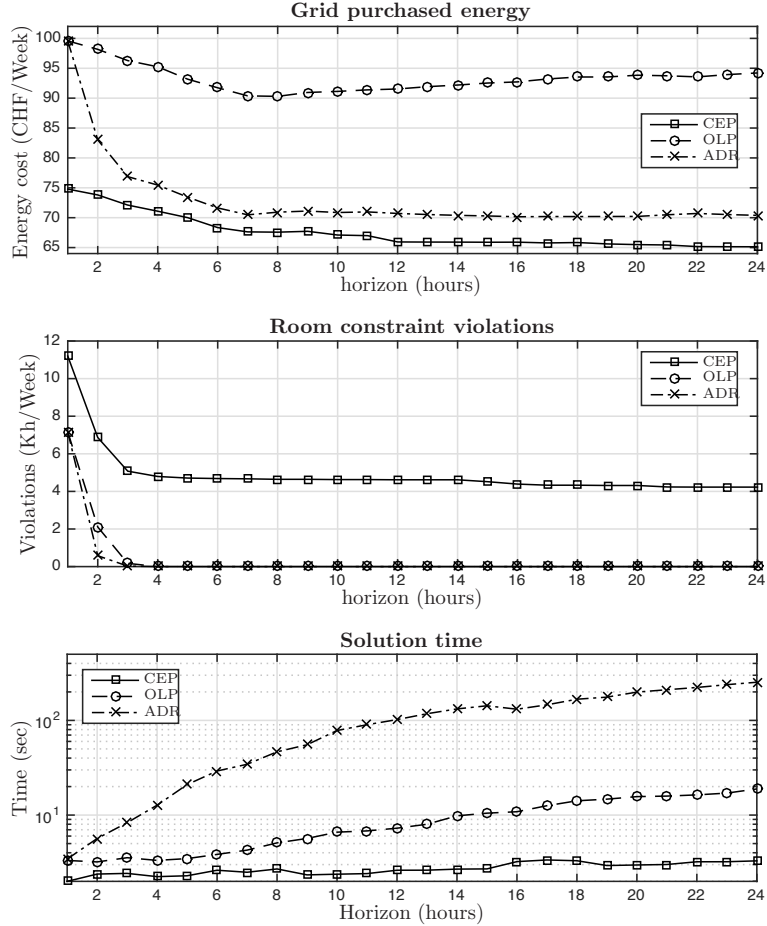}
		\caption{Profile trajectories for the mean room temperature of Building $ 1 $, and the total purchased grid energy by the community during a typical winter day.}
		\label{fig:tragectories}
	\end{center}
\end{figure}

To better visualize the behavior of the three approximations, the trajectories for the mean room temperature of Building $1$ (c.f. Table~\ref{tab:building_data}), and the amount of purchased grid energy, are depicted in Figure~\ref{fig:tragectories} for Monday, 1st January 2007. The CEP method operates very near comfort constraints and leads to frequent violations, while the OLP produces more conservative results by keeping the room temperature well inside the comfort range at the expense of consumption cost benefits. The poor performance of the CEP method with respect to the constraint violations can be explained when considering that the approximation assumes a deterministic evolution of the uncertain parameters. On the other hand, the adaptive nature of the ADR produces a good compromise between the optimistic decisions made by the CEP and the conservative decisions of the OLP. We remark that all three approximations utilize the load shifting capabilities of the battery by storing energy during the evening hours when electricity is cheaper, and deploying that energy in the early morning hours when the building needs to be brought back within the comfort range. As expected the energy produced from the photovoltaic unit is fully exploited to either reduce the grid energy consumption or charge the battery.

\subsection{Comparison of ADR and tuned-CEP methods}

We compare the stochastic ADR method with the deterministic \emph{tuned-CEP} one. We denote by {tuned-CEP}, the CEP method in which the comfort bounds of the building, $ i\in \mc B $, at time $ t \in \mc T $, have been tightened by a fixed constant $ c_b \ge 0 $ Celcius degrees, as follows,
\begin{equation*}
	\left.\begin{array}{r@{\,}l}
		\min\{\bs{x}_{t,i} - \text{lb}_{t,i},\; \text{ub}_{t,i} - \bs{x}_{t,i}\} &\geq c_b.
	\end{array}\right.
\end{equation*}
This constant $ c_b $ is computed as the minimum constraint tightening value for which the ADR and CEP methods achieve the same amount of comfort bounds violations over the simulation horizon. The {tuned-CEP} can be seen as a practical ad-hoc solution to deal with the increased constraint violations of the CEP method. 
\begin{table}[ht]
	\centering
	\renewcommand{\arraystretch}{1.0}
	\begin{tabular}{c|ccc}
		\multicolumn{1}{c}{~}&\multicolumn{3}{c}{\emph{\textbf{Winter}}} \\
		\hline
		\multicolumn{1}{c|}{Method} & \multicolumn{1}{c}{Cost (p.u.)} & \multicolumn{1}{c}{Tightening $ c_b $ (\degree C)} & Basis  (CHF/Week)\\
		\hline 
		\hline 
		CEP & $(1.00,\; 0.21)$ & [0.38,\,0.56]  & 327.56 \\
		ADR & $(0.94,\; 0.19)$ & 0  & \\
		\hline
	\end{tabular}
	\vspace*{2mm}
	\caption{Comparison of ADR and tuned-CEP methods. }
	\label{tab:resultsComp}
\end{table}

We conducted receding horizon simulations for $12$ consecutive weeks of the winter $2007$, and we report in \mbox{Table \ref{tab:resultsComp}} the empirical mean and standard deviation of the recorded weekly costs. For each of these weeks, we computed by repetitive trials the respective constraint tightening constant, $ c_b $, which ranges from \mbox{$ 0.38 $\degree C} and \mbox{$ 0.56 $\degree C}. The fact that $ c_b $ varies considerably from week to week, also suggests that tuning the CEP method is not a straightforward procedure. Additionally, we observe that the cost of purchased energy associated with the tuned-CEP solution method is considerably higher than the ADR one. Therefore, the inherent ability of the ADR method to deal with the uncertainty allows it to achieve an efficient trade-off between constraint violations and energy cost without the extra effort of tuning.

\subsection{District energy benefit}

Finally, we investigate the potential energy cost savings that can be obtained by considering diverse, heterogeneous buildings in a district. We define the heterogeneity on buildings by means of operational and construction characteristics. Specifically, we distinguish two operation types for the buildings, commercial (COM) and residential (RES). The day-night room temperature bounds defining these types are summarized in Table~\ref{tab::varyingBounds}. We select building (Bd) 1 and 3 of Table~\ref{tab:building_data} to be classified as commercial, while building 2 and 4 as residential. Notice that these buildings also differ in their construction characteristics.
\begin{table}[ht!]
	\renewcommand{\arraystretch}{1.0}
	\centering
	\begin{tabular}{|c|cc|cc|}
		\multicolumn{1}{c}{{}} & \multicolumn{2}{|c|}{COM Bounds} & \multicolumn{2}{c|}{RES Bounds} \\\hline 
		Time & Lower & Upper &  Lower & Upper \\\hline
		06:00 - 09:00  & 15\degree C & 30\degree C & 21\degree C & 25\degree C  \\
		09:00 - 19:00  & 21\degree C & 25\degree C & 15\degree C & 30\degree C \\
		19:00 - 23:00  & 15\degree C & 30\degree C & 21\degree C & 25\degree C  \\
		23:00 - 06:00  & 15\degree C & 30\degree C & 15\degree C & 30\degree C  \\
		\hline
	\end{tabular}
	\vspace*{1ex}
	\caption{Comfort constraints bounds for commercial and residential buildings. }
	\label{tab::varyingBounds}
\end{table}

We conduct five receding horizon simulations during the first week of January 2007. In each of these experiments, we consider a different combination of residential and commercial buildings from Table~\ref{tab:building_data}, and we report the cost of purchased grid energy for the CEP, ADR and OLP methods in Table \ref{tab::districtEnergyBen}. Notice that in each simulation experiment, we appropriately scale the hub devices with respect to the number of buildings that are connected to the energy hub. 
\begin{table}[ht!]
	\centering
	\renewcommand{\arraystretch}{1.0}
	\begin{tabular}{r|ccc|c}
		{} & ADR & OLP & CEP & Basis  (CHF/Week) \\
		\hline
		RES (Bd 2) 			& 1.13 & 1.39 & 1.00 & 56.37 \\
		RES (Bd 4) 			& 1.15 & 1.41 & 1.02 & \\
		COM (Bd 1) 			& 1.28 & 1.47 & 1.23 & \\
		COM (Bd 3) 			& 1.31 & 1.49 & 1.26 & \\
		RES+RES (Bd 2+4) 	& 2.27 & 2.80 & 2.02 & \\
		COM+COM (Bd 1+3)	& 2.58 & 2.96 & 2.49 & \\
		COM+RES (Bd 1+2)	& \textbf{1.51} & \textbf{1.90} & \textbf{1.37}&\\\hline
	\end{tabular}
	\vspace*{1ex}
	\caption{Cost of purchased energy for various district configurations during a typical winter week.}
	\label{tab::districtEnergyBen}
\end{table} 

We observe that the cost benefits obtained by merging buildings with different construction characteristics but similar operation plans can be limited. On the contrary, significant gains are obtained by aggregating buildings with dissimilar operation plans. There are two reasons for this: $ (i) $ The energy shifting mainly occurs among buildings with complementary demand profiles. In this occasion, the efficient, but of limited capacity, devices of the hub (e.g., heat pump) are fully utilized during the course of day; $ (ii) $ The energy produced by the photovoltaics units is better exploited through the storage when the peak demands of buildings do not coincide over time. This flexibility on the demand profiles results to a considerably less purchased energy from the grid.

\section{Conclusion}

This paper presents a unified data-driven control framework for the problem of cooperatively managing the aggregated energy demands of buildings in a district. It exploits the available historical data to train linear models with additive uncertainties that effectively capture the evolution of the stochastic processes in the system. The underlying distributions of these additive uncertainties are shown to belong to Gaussian families of distributions which are constructed off-line using the historical data information. We exploit the simple structure of these sets to approximate the resulting robust stochastic optimization problem by a finite dimensional linear program. This program is {tractable} and its complexity scales polynomially. This is particularly important when addressing large scale problems such as cooperative building energy management. An extensive simulation study based on realistic data demonstrates the efficacy of the proposed method. Among other things, our study demonstrates that additional cost benefits can be obtained by aggregating buildings with heterogeneity in the demand profiles.

As future work, we note that our problem has a decoupled structure. Indeed, the linear structure of the objective and and the weakly coupled structure of the constraints can be exploited  in a number of distributed and decentralized optimization algorithms. Additionally, the problem's decoupled structure can be exploited by recently developed optimization algorithms such as the alternating direction method of multipliers (ADMM), for a fast numerical solution of the linear optimization problem.

\section*{Acknowledgments}

The authors would like to thank the Building Science and Technology Laboratory (EMPA) for providing the building characteristics and the occupancy data, and MeteoSwiss for making available the weather forecasts and realizations. The authors would also like to thank Annika Eichler, Marc Hohmann and Ben Flamm for fruitful discussions on the topic.

\section*{Appendix}
\setcounter{equation}{0}
\renewcommand{\theequation}{\Alph{subsection}.\arabic{equation}}
\subsection*{Proof of Proposition \ref{prop::robGaus}}\label{app::1}

We construct the family of distributions, $ \mc P_{t} $ as the intersection of $ (i) $ the family, $ \mc P^\chi_{t} $, of Gaussian distributions with unknown mean and bounds on variance, and $ (ii) $ the family, $ \mc P^{st}_{t} $, of Gaussian distributions with unknown variance and bounds on the mean, as follows:
\begin{equation*}
	\begin{array}{r}
		\mc{P}_{t} =  \{\mb{P}_{t} \textrm{ such that } (\mb{P}_{t}\in\mc{P}^{\chi}_{t}) \wedge (\mb{P}_{t} \in \mc{P}^{\mt{st}}_{t})  \},
	\end{array}                     
\end{equation*}
where,
\begin{equation*}
	\begin{array}{r}
		\mathcal{P}^{\chi}_{t} = \{ \mathbb{P}_{t} \textrm{ is Gaussian and condition \eqref{eq::uncert::varianceBounds} is satisfied} \},\\[1ex]
		\mathcal{P}^{st}_{t} = \{ \mathbb{P}_{t} \textrm{ is Gaussian and condition \eqref{eq::uncert::meanBounds} is satisfied} \}.
	\end{array}
\end{equation*}
Given the data, $ \mc S_{t} $, equations \eqref{eq::uncert::varianceBounds} and \eqref{eq::uncert::meanBounds} also provide the confidence levels associated with the respective ambiguity sets,
\begin{equation*}
	\begin{array}{r}
		\bs{\mt{P}}_{\mc S_{t}}(\mathbb{P}_{t} \in \mathcal{P}^{{\chi}}_{t}) \ge 1-\delta^{\chi}_{t},\\[1ex]
		\bs{\mt{P}}_{\mc S_{t}}(\mathbb{P}_{t} \in \mathcal{P}^{{st}}_{t}|\,\mathbb{P}_{t} \in \mathcal{P}^{{\chi}}_{t}) \ge 1-\delta^{{st}}_{t}\,.
	\end{array}
\end{equation*}
Notice that the bounds in \eqref{eq::uncert::meanBounds} are derived under the assumption that the sample variance is following the chi-square distribution. Therefore, the confidence level associated with the construction of the distribution family, $ \mc P_t^{st} $, is given by the conditional probability. 

We use Bayes rule to evaluate the significance level, $ \delta_{t} $, associated with the ambiguity set $ \mc P_{t} $, as follows:
\begin{equation*}
	\begin{array}{r@{\,}l}
		\bs{\mt{P}}_{\mc S_{t}}(\mathbb{P}_{t} \in \mathcal{P}_{t}) =& \bs{\mt{P}}_{\mc S_{t}}(\mathbb{P}_{t} \in  \mathcal{P}^{\chi}_{t} \wedge \mathbb{P}_{t} \in \mathcal{P}^{\mt{st}}_{t})\\[1ex]
		=& \bs{\mt{P}}_{\mc S_{t}}(\mathbb{P}_{t} \in \mathcal{P}^{\mt{st}}_{t}\,|\, \mathbb{P}_{t} \in \mathcal{P}^{\chi}_{t} ) \cdot \bs{\mt{P}}_{\mc S_{t}}(\mathbb{P}_{t} \in \mathcal{P}^{\chi}_{t})\\[1ex]
		\ge & (1-\delta^{\mt{st}}_{t})(1-\delta^{\chi}_{t}) \\[1ex]
		= & 1 - \delta^{\mt{st}}_{t} -\delta^{\chi}_{t} + \delta^{\mt{st}}_{t} \delta^{\chi}_{t} \\[1ex]
		\ge & 1 - \delta^{\mt{st}}_{t} -\delta^{\chi}_{t} \\[1ex]
		= & 1 - \delta_{t}\,.
	\end{array}
\end{equation*}
with the significance level, $ \delta_t $, given as, $ \delta_t =\delta^{\mt{st}}_{t} +\delta^{\chi}_{t}. \hfill\blacksquare $

\subsection*{Proof of Proposition \ref{prop::strCon}}\label{app::3}
The uncertain vector $ \bs w $ comprises $ |\mc T| |\mc D| $ uncorrelated components. Therefore, we seek a bounded set $ \widehat{W} $ that is composed by $ 2 |\mc T| |\mc D| $ inequalities that upper and lower bound every component $ i =1,\ldots, |\mc T| |\mc D| $, of $ \bs w $. Denoting by $ \overline{w}_{t,i} $ and $ \underline{w}_{t,i} $ the upper and lower bound, respectively, the set $ \widehat{W} $ can be written as,
\begin{equation*}
	\widehat{W} = \big\{\bs w \in \mb{R}^{|\mc T| \,|\mc D|} \text{ : } \underline{w}_{t,i} \le w_{t,i} \le \overline{w}_{t,i},\,\forall t\in \mc T,\,i\in \mc D\big\}.
\end{equation*}
We require that,
\begin{equation*}
	\infimum_{\mb{P}  \in \mc{P} } \mb{P}(\bs w \in \widehat{W}) \ge 1-\epsilon\,.
\end{equation*}
We exploit the conservative Bonferonni approximation \cite{Nemirovski2006} to decouple the joint chance constraints into a set of individual ones, as follows:
\begin{equation*}\label{eq::bonfAppr}
	\begin{array}{r@{\,}l}
		\mb{P} (\bs w \notin \widehat{W}) &= \displaystyle\sum\limits_{i \in \mc D} \sum \limits_{t \in \mc T} \left(\mb{P} (w_{t,i}\ge \overline{w}_{t,i}) + \mb{P} (w_{t,i}\le \underline{w}_{t,i})\right)\\
		&=\displaystyle\sum\limits_{i \in \mc D} \sum \limits_{t \in \mc T} \left(\overline{\beta}_{t,i}+\underline{\beta}_{t,i}\right) = \epsilon,
	\end{array}
\end{equation*}
with $ \overline{\beta}_{t,i},\,\underline{\beta}_{t,i} \ge 0 $. Notice that in this case the Bonferonni approximation is exact since we are considering mutually exclusive events (i.e. $ \mb{P} (w_{t,i}\ge \overline{w}_{t,i} \wedge w_{t,i}\le \underline{w}_{t,i}) = 0 $).

\noindent Assume that $ w_{t,i} \sim \mc N(\mu_{t,i},\sigma^2_{t,i}) $, then the chance constraints,
\begin{equation*}
	\begin{array}{l}
		\mb{P} (w_{t,i}\le \overline{w}_{t,i}) \ge 1- \overline{\beta}_{t,i},\\[1ex]
		\mb{P} (w_{t,i}\ge \underline{w}_{t,i}) \ge 1- \underline{\beta}_{t,i},
	\end{array}
\end{equation*}
are equivalently reformulated (see \cite{LaGhBhJo03}), as follows,
\begin{equation*}
	\left. \begin{array}{@{}r@{}l}
		& {w}_i \ge { \mu}_i - \Phi^{-1}\big(1-\underline{\beta}_{t,i}\big) \, \cdot{\sigma}_{i},\\[1ex]
		&{w}_i \le { \mu}_i + \Phi^{-1}\big(1-\overline{\beta}_{t,i}\big) \,\cdot {\sigma}_{i}. 
	\end{array}\right.
\end{equation*}
In this setting, the robustified individual chance constraints,
\begin{equation*}
	\begin{array}{@{}l}
		\inf\limits_{\mb P \in \mc P}\mb{P} (w_{t,i}\le \overline{w}_{t,i}) \ge 1- \overline{\beta}_{t,i},\\[1ex]
		\inf\limits_{\mb P \in \mc P}\mb{P} (w_{t,i}\ge \underline{w}_{t,i}) \ge 1- \underline{\beta}_{t,i},
	\end{array}
\end{equation*}
are equivalently be reformulated as follows,
\begin{equation*}
	\left. \begin{array}{r@{\,}l}
		&\displaystyle {w}_{t,i} \ge { \mu}_{t,i} - \Phi^{-1}\big(1-\underline{\beta}_{t,i}\big) \, \cdot{\sigma}_{t,i}\,\\[1ex]
		&\displaystyle{w}_{t,i} \le { \mu}_{t,i} + \Phi^{-1}\big(1-\overline{\beta}_{t,i}\big) \,\cdot {\sigma}_{t,i} 
	\end{array}\right \rbrace \begin{array}{l}
	\forall \underline{\mu}_{t,i} \leq \mu_{t,i} \leq \overline{\mu}_{t,i},\\[1ex]\forall \underline{\sigma}_{t,i} \leq \sigma_{t,i} \leq \overline{\sigma}_{t,i}.
\end{array}
\end{equation*}
Exploiting the linearity of these constraints, we can explicitly compute the upper and lower bounds, as follows,
\begin{equation*}
	\left. \begin{array}{r@{\,}lr}
		&\displaystyle \underline{w}_{t,i} = \underline{ \mu}_{t,i} - \Phi^{-1}\big(1-\underline{\beta}_{t,i}\big)\,\cdot \overline{\sigma}_{t,i},&\\[1ex]
		&\displaystyle \overline{w}_{t,i} = \overline{ \mu}_{t,i} + \Phi^{-1}\big(1-\overline{\beta}_{t,i}\big) \, \cdot \overline{\sigma}_{t,i}.&
	\end{array}\right. 
\end{equation*}
This analysis provides the constructive proof for the structure of $ \widehat{W} $ in \eqref{eq::setW} $ \hfill \blacksquare $

\subsection*{Dynamics and constraints of the energy hub devices}\label{app::enHub}
We choose to present the devices characteristics sized per building (Bd). In particular, we model the chiller, boiler and heat pump using a coefficient of performance \cite{evins2014new}, which gives rise to the following constraints:
\begin{equation*}
	\begin{array}{l}
		\bs{u}^{\text{out}}_{t,\text{chiller}} = 0.7\bs{u}^{\text{in}}_{t,\text{chiller}},\\
		\bs{u}^{\text{out}}_{t,\text{boiler}} = 0.9\bs{u}^{\text{in}}_{t,\text{boiler}},\\
		\bs{u}^{\text{out}}_{t,\text{HP}} = 3\bs{u}^{\text{in}}_{t,\text{HP}},\\
	\end{array}
\end{equation*}
where $ \bs u_{t,i} = [\bs u^{\text{in}}_{t,i},\bs u^{\text{out}}_{t,i}] $
with \mbox{$\bs{u}^{\text{in}}_{t,i},\bs{u}^{\text{out}}_{t,i}$}, denoting the input and output power flows on the $ i $-th device, respectively. The capacities of the conversion units are given as follows,
\begin{equation*}
	\begin{array}{l}
		0\le \bs{u}^{\text{out}}_{t,\text{chiller}} \le 20\,\text{kW/Bd},\\
		0 \le \bs{u}^{\text{out}}_{t,\text{boiler}} \le 25\,\text{kW/Bd}\\
		0 \le \bs{u}^{\text{out}}_{t,\text{HP}} \le 5\,\text{kW/Bd},\\
	\end{array}
\end{equation*}
We consider a South oriented photovoltaic array with maximum output of $4.10\,\text{kW/Bd}$. The photovoltaic dynamics are generated by linearizing the non-linear model of \cite{fakham2007control} with respect to the ambient temperature and the solar radiation, 
\begin{equation*}
	0 \leq \bs{u}^{\text{out}}_{t,\text{PV}} \leq 0.1280  -0.0019\bs{\xi}_{t,\text{AT}} + 3.7\bs{\xi}_{t,\text{SRS}},
\end{equation*}
for all $t\in\mc T$. Note that the units of  $\bs{\xi}_{t,\text{SRS}}$ are measured in $\text{kW}/\text{m}^2$ and can typically take values $\bs{\xi}_{t,\text{SRS}}\in[0,1]$.

We consider a  lead-acid battery \cite{vrettos2011operating}, with a $5\text{kW/Bd}$ capacity, giving rise to the following linear dynamical system:  
\begin{equation*}
	\begin{array}{@{}r@{\,}l}
		\bs{x}_{t+1} = 
		\begin{pmatrix}
			0.51 & 0.22\\
			0.47 & 0.78
		\end{pmatrix}
		\bs{x}_{t} +
		\begin{pmatrix}
			0.61\\
			0.25
		\end{pmatrix}
		\bs{u}^{\text{in}}_{t}
		+  
		\begin{pmatrix}
			-0.83\\
			-0.39
		\end{pmatrix}
		\bs{u}^{\text{out}}_{t},
	\end{array}
\end{equation*}
where the states and control are constrained by:
\begin{equation*}\label{battery}
	\left.
	\begin{array}{l}
		0 \leq \bs{u}^{\text{in}}_{t},\,\bs{u}^{\text{out}}_{t} \leq 8,\\
		1\leq \bs{x}_{t,1} + \bs{x}_{t,2}\leq 5,\,\bs{x}_{t} \geq 0 ,\\
		0.62\bs{x}_{t,1} + 0.27\bs{x}_{t,2} - \bs{u}^{\text{out}}_{t} \geq 0, \\
		0.84\bs{x}_{t,1} + 0.37\bs{x}_{t,2} + \bs{u}^{\text{in}}_{t} \leq 2.58,\\
		0.73\bs{x}_{t,1} + 0.73\bs{x}_{t,2} + \bs{u}^{\text{in}}_{t} \leq 3.66.\\
	\end{array}
	\right.
\end{equation*}

To this end, we provide the equations describing the interconnection of the energy hub with the building community. In particular, the electricity balancing constraint is given by,
\begin{equation*}
	\bs{p}_{t,\text{elect}} + \sum_{i\in\text{E}_+} \bs{u}^{\text{out}}_{t,i} = \sum_{i\in\text{E}_-} \bs{u}^{\text{in}}_{t,i} + \bs{d}_{t,\text{elect}},
\end{equation*}
where $\text{E}_+ = \{$Photovoltaics, Battery$\}$ and $\text{E}_- = \{$Heat Pump, Chiller, Boiler, Battery$\}$. Similarly, the heating and cooling energy balancing constraints are given by,
\begin{equation*}
	\sum_{i\in\text{H}_+} \bs{u}^{\text{out}}_{t,i} = \bs{d}_{t,\text{heat}},\text{ and }\bs{u}^{\text{out}}_{t,\text{chiller}} = \bs{d}_{t,\text{cool}},
\end{equation*}
where $\text{H}_\text{+} = \{$Heat Pump, Boiler$\}$. Finally, the demand for electricity, heating and cooling of the building community is given by
\begin{equation*}
	\begin{array}{@{}lll}
		\bs{d}_{t,\text{elect}} = \displaystyle \sum_{i \in \mc B} \bs{u}_{t,i,\text{AHU}}, \\
		\bs{d}_{t,\text{cool}} = \displaystyle \sum_{i \in \mc B}\bs{u}_{t,i,\text{TABS}},\\
		\bs{d}_{t,\text{heat}} = \displaystyle \sum_{i \in \mc B}\big( \bs{u}_{t,i,\text{radiator}} + \bs{u}_{t,i,\text{TABS}}\big),
	\end{array}
\end{equation*}
where TABS are typically used both for heating and cooling of the buildings.

	\bibliographystyle{unsrt}
	\bibliography{darivianos_abrv,Papers,darivianos}
\end{document}